\documentclass{amsart}
\usepackage{amsmath,amssymb,amscd,amsthm}
\usepackage[matrix,arrow]{xy}

\newcommand{\N}         {\mathbb N}
\newcommand{\Z}         {\mathbb Z}

\newcommand{\R}         {\mathbb R}
\newcommand{\C}         {\mathbb C}

\newcommand{\F}         {\mathcal{F}}

\newcommand{\CP}        {\mathbb{CP}}

\newcommand{\Spinc}     {\mathrm{Spin}_{\C}}
\newcommand{\ra}        {\rightarrow}

\newcommand{\tensor}    {\otimes}

\newcommand{\union}{\cup}

\newcommand{\iso}       {\cong}

\newcommand{\Mhat}      {\hat{M}}
\newcommand{\Mred}      {M_{\mathrm{red}}}
\newcommand{\hso}       {h_{S^{1}}}
\newcommand{\Hso}       {H_{S^{1}}}
\newcommand{\Kso}       {K_{S^{1}}}
\newcommand{\ee}        {\tilde{e}}

\DeclareMathOperator{\Ker}{Ker}
\DeclareMathOperator{\Max}{Max}
\DeclareMathOperator{\Min}{Min}

\newtheorem{thm}{Theorem}
\newtheorem{Def}{Definition}
\newtheorem{proposition}{Proposition}

\newtheorem{lemma}{Lemma}

\begin{document}

\title{Cohomological Localization for Manifolds with Boundary}

\author{David S. Metzler}

\begin{abstract}
For an $S^{1}$-manifold with
boundary, we prove a localization formula applying
to any equivariant cohomology theory satisfying a certain
algebraic condition. We show how the localization result of Kalkman
and a case of the quantization commutes with reduction theorem
follow easily from the localization formula.
\end{abstract}

\maketitle


\section{Introduction}
\label{sec:intro}

In the paper \cite{Kalkman:1995Coho} on cohomology rings of symplectic
quotients,
Kalkman uses the Cartan model of ordinary (Borel) equivariant cohomology
to prove a localization result for $S^{1}$-manifolds with boundary.
In this paper we extend Kalkman's result 
to generalized equivariant cohomology theories satisfying
a certain algebraic condition.
In particular, the result applies to equivariant $K$-theory, and
in that setting it generalizes the version of ``quantization
commutes with reduction'' proved in \cite{Metzler:1996Kcut}.

We describe the basic goal briefly. 
Let $M$ be a compact oriented manifold with nonempty boundary, 
and with a circle action that is free on $\partial M$. Denote the inclusion
by $i: \partial M \ra M$. Denote the set of fixed point components by $\F$.
Let $X = \partial M/S^{1}$ and let $p: \partial M \ra X$ be the quotient map.
Note that this is a principal circle bundle. 

Let $h_{S^{1}}(-)$ be a (generalized) multiplicative, complex oriented, 
$S^{1}$-equivariant cohomology theory,
with corresponding\footnote{See section \ref{sec:localmhat}.} 
non-equivariant cohomology theory $h(-)$. Denote pullback
maps in cohomology by $f^{*}$ and Gysin maps (push-forwards) by $f_{!}$.
Assume that $X$ is $h$-oriented, so
that the pushforward $(\pi_{X})_{!}$ is defined,
where $\pi_{X}$ is the unique map $X \ra *$.
Our goal is a fixed point formula for the composed map
\begin{equation}\label{divideintegrate}
  \xymatrix{ {\kappa: \hso(M)} \ar[r]^-{i^{*}} 
     & \hso(\partial M) \ar[r]^-{(p^{*})^{-1}} 
     & h(X) \ar[r]^-{(\pi_{X})_{!}}
     & h(*).
  }
\end{equation}
This is given in Theorem \ref{mainthm}.

In Section \ref{sec:geometric} we present a simple
geometric construction that converts the boundary into a
fixed point component. In Section \ref{sec:localmhat} 
we apply the usual localization theorem to the new manifold.
In Section \ref{sec:residuemaps} we introduce the idea of a
residue map, which requires an additional assumption on the 
coefficient ring of the equivariant cohomology theory in question, 
and use the residue in Section \ref{sec:localizationformula} 
to get the desired formula. 
In Section \ref{sec:firstappl} we first show how this formula
specializes to Kalkman's in the case of Borel cohomology.
We then show that in $K$-theory it leads to a proof of a
restricted case of the ``quantization commutes with reduction''
result in symplectic geometry.

In a forthcoming paper 
we will investigate Theorem \ref{mainthm} in
other cohomology theories such as cobordism and elliptic cohomology.
We will also give a similar theorem for actions of $SU(2)$.

\section{A Geometric Construction}
\label{sec:geometric}

Given an $S^{1}$-manifold $M$, with notation as in the introduction, 
we use a simple procedure analogous to blowing up
(and to symplectic cutting 
\cite{EL:1995SC}) to produce a new manifold without boundary
which adds $X$ to the fixed point set.

Note: throughout the paper, by $S^{1}$ we really mean $S^{1} \subset \C$,
with the attendant explicit identification of $\mathrm{Lie}(S^{1})$
with $\R$ ($\exp(t) = e^{it}, \: t \in \R$). This is important for
orientation issues.

Since the principal stratum for a compact group action is open and dense,
there is a neighborhood $U$ of $\partial M$ on which the action of
$S^{1}$ is free. In fact, we can assume that $U$ is invariant and
that there is an equivariant homeomorphism 
$\phi: U \iso \partial M \times [0,\infty)$. 
Consider the action of $S^{1}$ on $\C$ of weight $-1$
(given $\lambda \in S^{1} \subset \C$, 
$\lambda \cdot z = \lambda^{-1}z$) and form the associated bundle
\[
  E = \partial M \times_{S^{1}} \C.
\]
(Here we divide by the diagonal action of $S^{1}$.)
This is an $S^{1}$-equivariant vector bundle, with $S^{1}$-action
defined by $\lambda \cdot [m,z] = [m,\lambda z]$.
Now $\phi$ induces a map 
\begin{align}\label{psi}
  \psi: U \setminus \partial M &\ra E \setminus X 
            \iso \partial M \times_{S^{1}} \C^{*} \\
        m &\mapsto [\phi_{1}(m),\phi_{2}(m)]         
\end{align}
which is easily seen to be a homeomorphism (since 
$\C^{*}/S^{1} \iso (0,\infty)$). It is also
$S^{1}$-equivariant, since
\begin{align*}
  \psi(\lambda \cdot m) 
       &= [\phi_{1} (\lambda \cdot m),
           \phi_{2} (\lambda \cdot m)] \\
       &= [\lambda \cdot \phi_{1} (m), \phi_{2} (m)] \\
       &= [\phi_{1} (m), \lambda \phi_{2} (m)]
\end{align*}
where the last equality comes from our taking the weight
$-1$ action to form the bundle $E$.

We use the map $\psi$ to glue, obtaining
\[
  \Mhat = (M \setminus \partial M) \union_{\psi} E.
\]
Clearly $\Mhat$ is a compact boundaryless $S^{1}$-manifold,
with fixed point set equal to the union of the fixed point
set of $M$ and $\partial M/S^{1} = X$.
$\Mhat$ is oriented, but we
have to be slightly careful about that. The orientation on
$M$ gives an orientation on $\partial M$ in the usual way.
Since $S^{1}$ is oriented, we get an orientation on the 
quotient $X$. 
We also have an orientation on the complex
vector bundle $E$ over $X$, but we need to use the
opposite orientation from the standard complex
orientation, in order to make $\psi$ orientation-preserving
Hence when we glue, we get an orientation
on $\Mhat$. 
Note that the normal bundle of $X$ in $\Mhat$ is exactly $E$,
oriented by the conjugate complex orientation.

In the simplest example, $M$ is the upper hemisphere of $S^{2}$ with 
the standard rotational action of $S^{1}$. Then $\Mhat$ is clearly $S^{2}$ with
the analogous action. We note that, given any orientation of
$S^{2}$, the circle action rotates in the positive direction
(with respect to the orientation) at one pole, and in the
negative direction at the other. This is an example of
the orientation issue discussed above.

The above description is the best for seeing the geometry of $\Mhat$
(in particular that it is a smooth manifold),
but an alternate description shows that it is a very natural topological
construction as well. 
\begin{proposition}\label{prop:hotpushout}
$\Mhat$ is homeomorphic to the homotopy pushout of the diagram
\[
  \xymatrix{ \partial M \ar@{^{(}->}[r]^{i} \ar[d]_{p} & M \\
           X.  
  }
\]
\end{proposition}
\begin{proof}
To take the homotopy pushout we must replace $p$ by its mapping cylinder  
\[
  \xymatrix{ \partial M \ar@{^{(}->}[r]^-{j} 
    & \frac{(\partial M \times I) \union X}
         {(m,1) \sim p(m)} = \mathrm{Cyl}(p) \iso D(E)
  } 
\]
where the last term denotes the disc bundle of $E$.
Hence when we take the ordinary pushout 
\[
 \frac{M \union \mathrm{Cyl}(p)}{i(m) \sim j(m)} 
\]
we are simply gluing $M$ to the disc bundle of $E$, which gives a
space homeomorphic to $\Mhat$.
\end{proof}

\section{Localization on $\Mhat$}
\label{sec:localmhat}

Let $\hso(-)$ be an $S^{1}$-equivariant, complex oriented 
multiplicative cohomology theory,
free and split over the cohomology theory $h(-)$, as in 
\cite{May:1996CBMS}. This implies the following facts:
the equivariant coefficient
ring $\hso = \hso(*)$ is an algebra over $h = h(*)$. 
Second,
for a trivial $S^{1}$-space $X$,
we have a natural isomorphism 
\begin{equation}\label{corresponding}
  \hso (X) \iso h(X) \tensor_{h} \hso.
\end{equation}
Third, for a free $S^{1}$-space $Y$, the projection map $p$ 
induces an isomorphism $h(Y/S^{1}) \iso \hso(Y)$; more precisely
the composed map
\[
  h(Y/S^{1}) \ra h(Y/S^{1}) \tensor \hso \iso \hso (Y/S^{1}) \ra \hso (Y)
\]
(the first map is just $\alpha \mapsto \alpha \tensor 1$),
which we will denote by $p^{*}$, is an isomorphism. Examples are Borel
cohomology \cite{Borel:1960TransSem}, equivariant 
$K$-theory \cite{Se:1968EK}, and (homotopical) equivariant bordism
\cite{May:1996CBMS}, \cite{Sinha:1999EquivBordism}.

We note a simple relation between $h_{S^{1}}(M)$ and
$h_{S^{1}}(\Mhat)$. 
\begin{proposition}\label{alphahatprop}
For every class $\alpha \in \hso(M)$ there is a unique
class $\hat{\alpha} \in \hso(\Mhat)$ such that
\begin{align}\label{alphahatproperties}
\hat{\alpha}|(M \setminus \partial M) &= \alpha|(M \setminus \partial M) \\
\hat{\alpha}|X &= (p^{*})^{-1} i^{*} \alpha.
\end{align}
\end{proposition}
\begin{proof}
Observe that the Mayer-Vietoris sequence
for $\Mhat$ which we get from Prop. \ref{prop:hotpushout} is
\begin{equation}\label{MVlong}
  \xymatrix{ {\cdots} \ar[r]^-{\partial} & \hso^{q} (\Mhat) \ar[r]
  & \hso^{q}(M) \oplus \hso^{q}(X) \ar[r] & \hso^{q}(\partial M) 
    \ar[r]^-{\partial} & {\cdots.}
  }
\end{equation}
However, since $\hso(X) \iso h(X) \tensor \hso$, and the map
$p^{*}: h(X) \ra \hso(\partial M)$ is an isomorphism, the map
\[
  \hso^{q}(M) \oplus \hso^{q}(X) \ra  \hso^{q}(\partial M) 
\]
in the Mayer-Vietoris sequence is surjective in every degree, forcing
the coboundary map to vanish.
Hence the long exact sequence splits into short exact sequences
\begin{equation}\label{MVshort}
  \xymatrix{ 0 \ar[r] & \hso^{q} (\Mhat) \ar[r]
  & \hso^{q}(M) \oplus \hso^{q}(X) \ar[r] & \hso^{q}(\partial M) 
    \ar[r] & 0.
  }
\end{equation}

Given $\alpha \in \hso(M)$, let 
$\beta = (p^{*})^{-1} i^{*} \alpha \in h(X)$.
Clearly the pair $(\alpha,\beta \tensor 1)$ maps to zero
in the Mayer-Vietoris sequence, and so defines a unique element
$\hat{\alpha} \in \hso(\Mhat)$, satisfying (\ref{alphahatproperties}).
\end{proof}
The class $\hat{\alpha}$ ``glues'' $\alpha$ to its ``quotient'' $\beta$. 

We now need to deal with $\hso$-orientations. In this case we
can define a reasonable notion of an $\hso$-orientation
on the manifold-with-boundary $M$. 
First recall that $\Mhat^{S^{1}} = M^{S^{1}} \cup X$, and that
the normal bundle of $X$ in $\Mhat$ is exactly $E$. 
We use its conjugate complex structure $\bar{E}$
to determine its $\hso$-orientation, to be consistent with the 
ordinary orientation, as discussed in Section \ref{sec:geometric}.
This leads to the 
\begin{Def}\label{orientationwithboundarydefinition}
  Let $\hso$ be an $S^{1}$-equivariant cohomology theory
  as above. Let $M$ be a manifold with boundary, with 
  $S^{1}$-action free on $\partial M$, and construct $\Mhat$ 
  as above. Assume that $X = \partial M/S^{1}$ has a given
  $h$-orientation and that each $F \in \F$ and each normal bundle 
  $\nu_{F}$ have given $\hso$-orientations.

  Then an $\hso$-\textbf{orientation} on $M$ compatible with this data
  is defined to be an orientation on $\Mhat$, compatible with
  the orientations on $F$ and $\nu_{F}$ in the usual way,
  and compatible with the orientation on $X \subset \Mhat$ 
  and the conjugate complex orientation on $E \iso \nu_{X}$.
\end{Def}

In the case of ordinary Borel cohomology, an $\hso$-orientation
on a manifold without boundary 
is simply an ordinary (geometric) orientation. From the discussion
in Section \ref{sec:geometric} it is clear that a geometric
orientation on $M$ gives a Borel $\hso$-orientation as defined above.

We seek a fixed point formula for the map
\begin{equation}\label{divideintegrateagain}
  \xymatrix{ \hso(M) \ar[r]^-{i^{*}} 
     & \hso(\partial M) \ar[r]^-{(p^{*})^{-1}} 
     & h(X) \ar[r]^-{(\pi_{X})_{!}}
     & h(*).
  }
\end{equation}

Any multiplicative equivariant cohomology theory has a localization 
theorem \cite{Dieck:1971localization}, \cite{KK:1991},
which gives a fixed point formula for 
\[
  \xymatrix{ {\hso(\Mhat)} \ar[r]^-{(\pi_{\Mhat})_{!}} 
    & {\hso.}
  }
\]
Specifically, it says that after sufficiently localizing (algebraically)
in the ring $\hso$ and the $\hso$-module $\hso(\Mhat)$, we have
\begin{equation}\label{localmhat}
(\pi_{\Mhat})_{!} \hat{\alpha} 
   = \sum_{F \in \hat{\F}} 
        (\pi_{F})_{!} \frac{i_{F}^{*} \hat{\alpha}}{\tilde{e}(\nu_{F})} 
\end{equation}
where $\hat{\F}$ is the set of fixed point components of $\Mhat$, 
$i_{F}: F \ra M$ is the inclusion map and $\tilde{e}$ is the $\hso$
Euler class (the restriction of the $\hso$ Thom class to the zero
section). The algebraic localization necessary 
inverts the following subset $S$ of
$\hso$. Let $\mathcal{R}$ be the set of finite dimensional,
$\hso$-oriented representations of $S^{1}$ with no trivial summand
(where a representation $V$ 
is regarded as an equivariant bundle over a point). Let
\[
  S = \{ \tilde{e}(V) \: | \: V \in \mathcal{R} \}.
\]
Then equation (\ref{localmhat}) holds in the ring $S^{-1} \hso$.

The $\hso$-orientation on $M$, and the resulting 
orientations on $\Mhat$, $X$, $F$, and $\nu_{F}$ allow
us to use the localization formula.
Note that the $S^{1}$-action on $\bar{E}$ has weight $-1$. 
We saw above that
$\hat{\alpha}|_{X} = (p^{*})^{-1} i^{*} \alpha$, hence
\begin{equation}\label{localmhatwithx}
(\pi_{\Mhat})_{!} \hat{\alpha} 
   = \sum_{F \in \F} 
        (\pi_{F})_{!} \frac{i_{F}^{*} \alpha}{\tilde{e}(\nu_{F})} 
     + (\pi_{X})_{!} \frac{(p^{*})^{-1} i^{*} \alpha} 
                          {\tilde{e}(\bar{E})}.
\end{equation}

This is our key formula, but note the presence of the unknown quantity
$(\pi_{\Mhat})_{!} \hat{\alpha}$. The next step is to get rid of it,
and get rid of 
the factor $(\tilde{e}(\bar{E}))^{-1}$. To do that we need a further
assumption about $\hso$, which goes beyond the usual axioms for
equivariant cohomology theories, but is satisfied, for example, 
by Borel cohomology and equivariant $K$-theory. 

\section{Residue Maps}
\label{sec:residuemaps}

First we need a notation for the universal Euler classes of the
theory. Let $\mathcal{L} \ra \CP^{N}$ ($N \in \N$) be the universal
line bundle, where both $\mathcal{L}$ and $\CP^{N}$ have
the trivial $S^{1}$-action. Let $\C_{m}$ be the weight $m$ irreducible
representation of $S^{1}$, given by $\lambda \cdot z = \lambda^{m}z$. Let 
$\mathcal{L}_{m} = \mathcal{L} \tensor \C_{m}$.
Denote the $\hso$-Euler class of $\mathcal{L}_{m}$ by
\[
  E_{m} = \tilde{e}(\mathcal{L}_{m}) 
                \in \hso (\CP^{N}) \iso \hso[y]/(y^{N+1}).
\]
Here $y = e(L)$.
(We suppress $N$ from the notation because it will be irrelevant
to what follows, as long as it is large enough. However we do not use 
$\CP^{\infty}$ because we need $y$ to be nilpotent.)

Denote the canonical localization map by $\varepsilon: \hso \ra S^{-1} \hso$.
We will need the existence of the following kind of map.
\begin{Def}\label{rhodef}
A \textbf{residue map} for $\hso(-)$ is an $h$-module homomorphism 
$\rho: S^{-1}\hso \ra h$ such that:
\begin{enumerate}
\item \label{rhokillsgood}
      $\rho \circ \varepsilon = 0$.
\item \label{rhoandeuler}
      The natural extension of $\rho$ to a map 
      $\rho: S^{-1}\hso[y]/(y^{N+1}) \ra h[y]/(y^{N+1})$ satisfies
      \[
         \rho ((E_{-1})^{-1}) = -1.
      \]
\end{enumerate}
\end{Def}

A simple naturality argument
shows that condition \ref{rhoandeuler} implies the following:
given a trivial $S^{1}$-space $X$ and 
$L$ an $S^{1}$-equivariant line bundle of weight $-1$ over $X$,
the map $1 \tensor \rho$ takes 
$(\tilde{e}(L))^{-1} \in h(X) \tensor S^{-1} \hso$ to 
$-1 \in h(X) \tensor h \iso h(X)$. This is how we will
use condition \ref{rhoandeuler}.

We now give the conditions under which such a residue map exists,
in terms of the the formal group law of $\hso$.
Let $F$ be the formal group law of $\hso$, i.e.
$F(X,Y) \in \hso[[X,Y]]$ is the formal power series such that
for two line bundles $L_{1},L_{2}$ on $X$,
\begin{equation}\label{formalgrouplaweulerproperty}
 \tilde{e}(L_{1} \tensor L_{2}) = F(\tilde{e}(L_{1}) \tensor \tilde{e}(L_{2})).
\end{equation}
It is well-known that 
\[
  F(X,Y) = \sum_{k,l=0}^{\infty} a_{k,l} X^{k} Y^{l} = X + Y + XY\hat{F}(X,Y)
\]
where $\hat{F}(X,Y) \in \hso[[X,Y]]$ and $\hat{F}(X,Y) = \hat{F}(Y,X)$.
Hence $a_{k0} = a_{0k} = \delta_{1k}$.
Let $e = \tilde{e}(\C_{-1})$. Then 
\begin{align}\label{eulercalcformalgrouplaw}
E_{-1} = \ee(\C_{-1} \tensor L) 
 &= F(e,y) \\  
 &= \sum_{k,l=0}^{\infty} a_{kl}e^{k}y^{l} \\
 &= \sum_{l=0}^{\infty} b_{l} y^{l}
\end{align}
where $b_{l} \in \hso$ are defined by the last equation, or
by directly calculating $\ee(\C_{-1} \tensor L) \in \hso(\CP^{N})$. 
Since $a_{k0} = a_{0k} = \delta_{1k}$, $b_{0}=e$.
Hence 
\begin{align}\label{inverteulerformal}
(E_{-1})^{-1} 
 &= e^{-1} [1 + \sum_{l=1}^{\infty} (b_{l}/e)y^{l}]^{-1} \\
 &= e^{-1} \sum_{k=0}^{\infty} 
             \left( - \sum_{l=1}^{\infty} (b_{l}/e)y^{l} \right)^{k} \\
 &= e^{-1} + \sum_{k=1}^{\infty} 
               -e^{-1} \left( \sum_{l=1}^{\infty} (b_{l}/e)y^{l} \right)^{k} \\
 &= e^{-1} + \sum_{k=1}^{\infty} c_{k}y^{k}
\end{align}
where the coefficients $c_{k} \in \hso[e^{-1}]$ 
are determined by the last equation. The sums are actually finite
for any given $N$.

Let $M$ be the $h$-submodule of $S^{-1}\hso$ generated by the
$c_{k}$, let 
\[
N = \frac{S^{-1}\hso}{\hso + M}, 
\]
and let $A$ be the $h$-submodule of $N$ generated by $e^{-1}$.
\begin{proposition}\label{existenceuniqueofresidue}
A residue map $\rho$ exists for $\hso(-)$ if and only if 
$e^{-1}$ is not torsion in $N$ as an $h$-module and the inclusion
$A \subset N$ splits as a map of $h$-modules.

If a residue map exists then the set of possible residue maps
is given by $\mathrm{Hom}_{h}(N/A,h)$.
\end{proposition}
\begin{proof}
 In the definition of a residue, clearly condition (\ref{rhokillsgood})  
 is satisfied iff $\rho$ defines a map from $S^{-1}\hso/\hso$ to $h$.
 Since 
  \[
   \rho((E_{-1})^{-1}) = \rho(e^{-1}) + \sum_{k=1}^{\infty} \rho(c_{k})y^{k}, 
  \]
 condition (\ref{rhoandeuler}) will be satisfied iff $\rho(c_{k}) = 0$ 
 for all $k \ge 1$ and $\rho(e^{-1}) = -1$. The former says that
 $\rho$ defines a map from $N = S^{-1}\hso/(\hso + M)$ to $h$,
 while the latter requires that this map split the inclusion of
 $A = \langle e^{-1} \rangle$ into $N$. The set of such splittings is just
 $\mathrm{Hom}_{h}(N/A,h)$. 
\end{proof}

In particular, when $h$ is a field, the conditions for existence
are automatically satisfied as long as $e^{-1} \not\in \hso + M$.

We now show that our two main examples have such a map $\rho$.
We will consistently work over $\C$, without further comment.

\subsection{Borel}\label{subsec:borel} 
Borel cohomology is defined by 
\[
  \Hso(M) = H(M \times_{S^{1}} ES^{1},\C).
\]
This is what is usually referred to in the geometry literature simply as
``equivariant cohomology.'' It can be computed by the deRham-theoretic 
Cartan model (see \cite{AB:1984MM}). We have
$\hso \iso \C[u]$, $e = -u$, $S = \{(mu)^{k} \: | \: m \ne 0, k \ge 0 \}$. 
Hence $S^{-1}\hso \iso \C[u,u^{-1}]$. The formal group law is the
additive one, $F(X,Y) = X + Y$, so $\hat{F} = 0$.
Therefore 
\begin{align}\label{invertforborel}
 (E_{-1})^{-1} 
   &= (-u + y)^{-1} \\
   &= -u^{-1} (1 - (y/u))^{-1} \\
   &= -u^{-1} + \sum_{k=1}^{\infty} y^{k}u^{-k-1}
\end{align}
so $c_{k} = u^{-k-1}$, $k \ge 1$, and 
$M = \langle u^{-k-1} \: | \: k \ge 1 \rangle$.
We have 
\[
h_{T} + M = \langle u^{l} \: | \: l \in \Z, l \ne -1 \rangle
\]
so $N \iso \langle u^{-1} \rangle$ and the inclusion $A \subset N$
is in fact an isomorphism. Hence there is a unique residue map
\[
 \rho \left( \sum_{k=-M}^{M'} a_{k}u^{k} \right) = a_{-1}, 
\]
or in complex-analytic terms, $\rho(f) = \mathrm{Res}_{u=0}(f(u) \: du)$.

\subsection{$K$-theory}\label{subsec:ktheory}
Here $h \iso \C$ and the equivariant coefficient ring is 
$R(S^{1}) \tensor \C \iso \C[u,u^{-1}] \iso \C[e,(1-e)^{-1}]$
where $e = \tilde{e}(\C_{-1}) = 1 - u^{-1}$. 
(See \cite{ASe:1968}, \cite{Se:1968EK}.)
$S$ is the multiplicative subset generated by 
$\{ 1-u^{k} \: | \: k \ne 0 \}$, so 
\[
S^{-1}\hso \iso \C[u,u^{-1},\{(u-\gamma)^{-1} \: | \: 
                               \gamma \text{~a root of unity~} \}].
\]
The formal group law is $F(X,Y) = X + Y - XY$, so $\hat{F}= -1$, or
$b_{0}=e$, $b_{1}=1-e$.
Therefore
\[
  (E_{-1})^{-1} = e^{-1} + \sum_{k=1}^{\infty} (-y)^{k} e^{-1}(e^{-1}-1)^{k}
\]
and
\begin{align}\label{ckforktheory}
  c_{k} = e^{-1}(1-e^{-1})^{k} 
   &= \frac{1}{1-u^{-1}}\left(1 - \frac{1}{1-u^{-1}} \right)^{k} \\
   &= \frac{-u}{(1-u)^{k+1}} \\ 
   &= \frac{1}{(1-u)^{k}} - \frac{1}{(1-u)^{k+1}}.  
\end{align}

Let $B$ be the vector subspace generated by 
$\{ (1-u)^{-k} \: | \: k \ge 1 \}$ and
let $D$ be generated by 
\[
\{ (u-\gamma)^{-k} \: | \: k \ge 1, 
    \gamma \text{~a root of unity}, \gamma \ne 1 \}.
\]
Let $\sigma: B \ra \C$ be defined by 
\[
  \sigma(\sum a_{k}(1-u)^{-k}) = \sum a_{k}.
\]
Then $S^{-1}\hso \iso \hso \oplus B \oplus D$, and
$M = \Ker \sigma$. Hence 
\[
  N = \frac{S^{-1}\hso}{\hso + M} \iso B/M \oplus D \iso \C \oplus D
\]
where the second isomorphism is induced by $\sigma$. 
Under this map 
\[
e^{-1} = \frac{1}{1-u^{-1}} = 1 - \frac{1}{1-u} \mapsto -1 \in \C \subset 
\C \oplus D.
\]
Hence any extension of $\sigma$ to $B \oplus D$ is a residue map.
The set of possible residue maps is naturally isomorphic to 
$D^{*}$.

There are two simple choices for the map $\rho$.
\begin{lemma}\label{rhoasresidueatunity}
 The map $\rho_{1}: \C(u) \ra \C$ defined by 
 $\rho(f) = - \mathrm{Res}_{u=1}(u^{-1}f(u) \: du)$ restricts to 
 $S^{-1}\hso$ to be the unique residue map vanishing on $D$.
\end{lemma}
\begin{proof}
 The map $\rho_{1}$ evidently vanishes on $\hso \iso \C[u,u^{-1}]$   
 and $D$ since every element of either set has no pole at $u=1$.
 On $B$,
\begin{align*}
  \frac{u^{-1}}{(1-u)^{k}} 
  &= (1-u)^{-k} [1-(1-u)]^{-1} \\ 
  &= (1-u)^{-k} \sum_{l=0}^{\infty} (1-u)^{l} \\ 
  &=  \sum_{l=0}^{\infty} (-1)^{l-k} (u-1)^{l-k}  
\end{align*}
so 
\[
\rho_{1}((1-u)^{-k}) = -\mathrm{Res}_{u=1}(u^{-1}(1-u)^{-k} \: du)
  = - (-1)^{-1} = 1 = \sigma((1-u)^{-k}).
\]
So $\rho_{1}$ extends $\sigma$ as desired, hence is the required
residue map.
\end{proof}

There is a more convenient residue map for the application
to quantization in Section \ref{sec:firstappl}. 
Let $L_{+} = \C[[u]][u^{-1}]$ and $L_{-} = \C[[u^{-1}]][u]$ be the
rings of positive and negative Laurent series. Given a
rational function $f \in \C(u)$, define $L_{+,0}(f)$
to be the constant term in the positive Laurent series for $f$,
and define  $L_{-,0}(f)$
to be the constant term in the negative Laurent series for $f$.
Define $\rho_{0,\infty}: \C(u) \ra \C$ by 
\[
  \rho_{0,\infty} = L_{+,0} - L_{-,0}.
\]
In complex analytic terms, 
\[
  \rho_{0,\infty}(f) 
    = (\mathrm{Res}_{u=0} + \mathrm{Res}_{u=\infty})(u^{-1} f(u) \: du)
\]

\begin{lemma}\label{rhoatzeroandinfinity}
 The map $\rho_{0,\infty}$ restricts to $S^{-1}\hso$
 to be a residue map, distinct from $\rho_{1}$.
\end{lemma}
\begin{proof}
  The residue theorem implies that for a function $f$ with poles only
  at $0,1,\infty$, $\rho_{1}(f) = \rho_{0,\infty}(f)$.
  Hence the two maps agree on $S^{-1}\hso/D$, so $\rho_{0,\infty}$
  is a residue map. However it is easy to see that they disagree on
  $D$. In fact a simple calculation shows that
\[
  \rho_{0,\infty}((u-\gamma)^{-k}) = (-\gamma)^{-k} \ne 0.
\]
\end{proof}

It is interesting that in both these cases, the equivariant
coefficient rings are integral domains, and the residue maps
extend easily to the entire field of fractions $\C(u)$.

\section{The Localization Formula}
\label{sec:localizationformula}

\begin{thm}\label{mainthm}
Let $\hso(-)$ be an $S^{1}$-equivariant cohomology theory, 
free and split over $h(-)$, and suppose that $\rho$ is a
residue map for $\hso(-)$. Let $M$ be an $S^{1}$-manifold
with boundary $\partial M$ on which $S^{1}$ acts freely,
and let $X = \partial M/S^{1}$. Assume that $M$ is $\hso$-oriented,
and define $\kappa$ as in the introduction. Then for every
$\alpha \in \hso(M)$,
\begin{equation}\label{mains1}
  \kappa(\alpha) 
      = \sum_{F \in \F} 
        (\pi_{F})_{!} \rho \left( \frac{i_{F}^{*} \alpha}
                                {\tilde{e}(\nu_{F})} \right).
\end{equation}
\end{thm}
\begin{proof}
Recall that ordinary localization on $\Mhat$ gave (\ref{localmhatwithx}):
\begin{equation*}
(\pi_{\Mhat})_{!} \hat{\alpha} 
   = \sum_{F \in \F} 
        (\pi_{F})_{!} \frac{i_{F}^{*} \alpha}{\tilde{e}(\nu_{F})} 
     + (\pi_{X})_{!} \frac{(p^{*})^{-1} i^{*} \alpha} 
                          {\tilde{e}(\bar{E})}.
\end{equation*}
Now $(\pi_{\Mhat}) \hat{\alpha} \in \hso$, so 
applying $\rho$ to (\ref{localmhatwithx}) gives (writing
$\rho$ in place of $1 \tensor \rho$)
\begin{align}\label{rhoapplied}
  0 &= \sum_{F \in \F} 
          \rho (\pi_{F})_{!} \frac{i_{F}^{*} \alpha}
                                  {\tilde{e}(\nu_{F})}
        + \rho (\pi_{X})_{!} \frac{(p^{*})^{-1} i^{*} \alpha}
                                  {\tilde{e}(E)} \\
    &= \sum_{F \in \F} 
          (\pi_{F})_{!} \rho \frac{i_{F}^{*} \alpha}
                                  {\tilde{e}(\nu_{F})}
        - (\pi_{X})_{!} (p^{*})^{-1} i^{*} \alpha
\end{align}
because $\rho$ is an $h$-module map.
The last term is just $\kappa(\alpha)$, so the theorem is proved.
\end{proof}

\section{First Applications}
\label{sec:firstappl}

First we show how this result agrees with and extends known
results in the Borel and $K$-theory cases. 

\subsection{Kalkman's Result}
\label{kalkmancase}
In this case the result is exactly Kalkman's formula, 
which he proves using the Cartan Model. We
compare it to the form given in \cite{GK:1996}.

We have already remarked on the fact that an ordinary geometric
orientation of $M$ gives exactly the kind of orientation we need
in the Borel theory.

Their map $r$ (the Kirwan map)
is defined to be $(p^{*})^{-1} \circ i^{*}$. Their
map ``res'' is exactly our $\rho$ for the Borel case.
Hence
\begin{align}\label{kalkmaneq}
  \int_{X} r (\alpha)  
      &= \kappa(\alpha) \\
      &= \sum_{F \in \F} 
        (\pi_{F})_{!} \rho \left( \frac{i_{F}^{*} \alpha}
                                {\tilde{e}(\nu_{F})} \right) \\
      &= \sum_{F \in \F} 
         \int_{F} \mathrm{res} \left( \frac{\alpha|_{F}}
                                {\tilde{e}(\nu_{F})} \right)
\end{align}
which is exactly Kalkman's localization formula.

\subsection{Quantization Commutes With Reduction}
\label{qrcase}
The most interesting application of the formula in $K$-theory is
in the case of a Hamiltonian circle action. Let $(M,\omega)$ be a
symplectic manifold with a Hamiltonian circle action, with moment map 
$\phi: M \ra \mathrm{Lie}(S^{1})^{*} \iso \R$.
Assume that the symplectic form is integral. Then
it is well-known that there is a ``prequantum'' line bundle
$L$ over $M$, with an action of $S^{1}$ lifting the given
one on $M$, whose Chern class is $[\omega]$.
Further, the action of $S^{1}$ on $L$ over a fixed point component
is simple to describe (see \cite{BK:1970}): on a component
$F$, the moment map takes a fixed, integral value $\phi(F)$.
The weight of the action of $S^{1}$ on the fibers of the line bundle
is exactly $-\phi(F)$. 

Any symplectic manifold has a natural isotopy class of complex
structures associated to it (\cite{McDuffSalamon:1995}). Choose such a 
structure. Since the circle action is symplectic, it preserves 
the complex structure. We will actually use the conjugate of this
complex structure as the $K_{S^{1}}$-orientation on $M$. 
(This is to agree with an index-theoretic definition of the
quantization---see \cite{Metzler:1996Kcut} for details.)
The \textit{quantization} of the space $M$ is defined to be the
pushforward $Q(M) := (\pi_{M})_{!}([L]) \in R(S^{1})$. 

When we have a Hamiltonian group action, it is natural to perform
symplectic reduction. The reduced space is defined to be 
\[
  M_{0} = \phi^{-1}(0)/S^{1}
\]
which is a smooth symplectic manifold if $0$ is a regular value of $\phi$
and the action of $S^{1}$ on the zero level set is free;
we will assume these conditions hold for the remainder of the
section. If $M$ is prequantizable then so is $M_{0}$,
where the prequantum line bundle on the reduced space is
obtained by restricting $L$ to the zero level set $\phi^{-1}(0)$, and then
identify this with a (non-equivariant) line bundle $L_{0}$ downstairs.

Since $M_{0}$ is symplectic, it has a complex structure, and
we again use the conjugate structure to give it a $K$-orientation.
The quantization of $M_{0}$ is the pushforward 
$(\pi_{M_{0}})_{!}(L_{0}) \in \Z$. 

\begin{thm}[\cite{DGMW:1995}]\label{QRthm}
Given a prequantizable Hamiltonian $S^{1}$ space
with a free action on the zero level set, we have
\[
  Q(M_{0}) = Q(M)_{0},
\]
where the right hand side stands for the multiplicity of the
trivial representation in the virtual representation $Q(M)$.
\end{thm}

We will derive this theorem from a more general result,
Theorem \ref{quantizationtheorem} below, which is an extension of a result
originally due to A. Canas da Silva, Y. Karshon and 
S. Tolman \cite{CanasKarshonTolman:1997QP}.
They consider three kinds of quantization of 
presymplectic manifolds, including the case of $\Spinc$-quantization.
We will phrase our generalization of their 
main result in $K$-theoretic terms.

The data for this version of quantization are: a compact $S^{1}$-manifold $M$, 
a fixed $\Kso$-orientation on $M$, and a class $\alpha \in \Kso(M)$.
The quantization of $M$ is then defined to be the pushforward 
\[
  Q(M) = (\pi_{M})_{!} \alpha \in \Kso(*) \iso R(S^{1}).
\]

To define a general version of reduction, we follow 
\cite{CanasKarshonTolman:1997QP} (with minor modifications in notation) 
and make the
\begin{Def}\label{reduciblehypersurface}
  A \textbf{reducible hypersurface} in $M$ is a co-oriented
  submanifold $Z$ of codimension one which is invariant under the
  $S^{1}$-action and on which this action is free. The
  \textbf{reduction} of $M$ at $Z$ is the quotient 
  $\Mred = Z/S^{1}$. $Z$ is \textbf{splitting} if 
  $M$ is the union of two (not necessarily connected)
  manifolds with boundary, say $M_{+}$, $M_{-}$, along their common boundary 
  $Z$, such that positive (negative) normal vectors
  to $Z$ point into $M_{+}$ ($M_{-}$). We then say that $Z$ \textbf{splits}
  $M$ into $M_{+}$, $M_{-}$.
\end{Def}

For example, if $\phi: M \ra \R$ is an invariant map, with $0$
a regular value, and if $S^{1}$ acts freely on $Z = \phi^{-1}(0)$, then
$Z$ splits $M$ into $M_{+} = \phi^{-1}([0,\infty))$
and $M_{-} = \phi^{-1}((-\infty,0])$. The main example we have
in mind is where $M$ is Hamiltonian and $\phi$ is a moment map.

Let $Z$ be a splitting, reducible hypersurface in $M$. We have exact
sequences 
\begin{equation}\label{reducedses}
  \xymatrix{ 0 \ar[r] & TZ \ar[r] & i^{*}TM \ar[r] & \nu_{Z} \ar[r] & 0 \\
             0 \ar[r] & V \ar[r] & TZ \ar[r] & \pi^{*}(T\Mred) 
                  \ar[r] & 0
  }
\end{equation}
where $V$, the vertical subbundle, is the same as the set of vectors that
are tangent to the orbits of $S^{1}$. Note that $V$ and $\nu_{Z}$
are both oriented trivial real line bundles. Hence
\[
  i^{*}TM \iso \pi^{*}T\Mred \oplus \R^{2}.
\]
The $\Kso$-orientation on $M$ hence gives a $\Kso$-orientation
on $\pi^{*}T\Mred$, hence a $K$-orientation on $\Mred$.

The class $\alpha$ descends to a class $\alpha_{\mathrm{red}} \in K(\Mred)$, 
first by restriction to
$Z$, and then by the isomorphism $\Kso(Z) \iso K(\Mred)$.
Hence we can define the quantization of $\Mred$ exactly
as above, 
\[
  Q(\Mred) = (\pi_{\Mred})_{!} \alpha \in \Z.
\]
This defines the quantization of the reduction of $M$; the 
reduction of the quantization of $M$ is defined to be 
$Q(M)_{0}$, the multiplicity of the trivial representation in
$Q(M)$. Theorem \ref{quantizationtheorem} will say that these
are equal, subject to certain conditions on $\alpha$ and the
normal bundles of the fixed point components of $M$.

In order to get the sharpest version of the theorem, we will
be a little tricky about orientations on the fixed point components.
Given a component $F \subset M^{S^{1}}$ with $F \subset M_{+}$,
equip the normal bundle $\nu_{F}$
with a complex structure by declaring all of the weights of the
$S^{1}$-action to be positive, and use the associated $\Kso$-orientation. 
Together with the given $\Kso$-orientation on $M$, 
this gives an $\Kso$-orientation on $F$. For $F \subset M_{-}$,
use the complex structure and associated $\Kso$-orientation
induced by declaring the weights to be negative.

It is not hard to see that the $\Kso$-orientation on $M$ gives 
a $\Kso$-orientation 
on $M_{+}$, $M_{-}$, in the sense of 
Definition \ref{orientationwithboundarydefinition}, 
compatible with the above orientations for the fixed point
components and the reduced space.

Given a connected component $F$ of $M^{S^{1}}$, let 
\[
 w(F) = \sum |\text{weights of $S^{1}$-action on $\nu_{F}$}|.
\]
Note that the absolute value of a weight is well-defined independent 
of the complex structure chosen. 

Recall that 
\[
  \Kso(F) \iso K(F) \tensor R(S^{1}) \iso K(F) \tensor \Z[u,u^{-1}].
\]
Hence for $\alpha \in \Kso(M)$, 
\[
i_{F}^{*} \alpha = \sum_{\Min(\alpha,F)}^{\Max(\alpha,F)} \alpha_{F,k} u^{n}
\]
for some integers $\Min(\alpha,F),\Max(\alpha,F)$. 
(Note that if $\alpha_{F}$ is homogeneous in $u$, then 
$\Min(\alpha,F)=\Max(\alpha,F)=\deg(i_{F}^{*}\alpha$).

We can now state the
\begin{thm}\label{quantizationtheorem}
  Let $M$ be a smooth, compact, $\Kso$-oriented $S^{1}$-manifold and 
  let $\alpha \in \Kso(M)$.
  Let $Z$ be a reducible hypersurface which splits $M$ into
  $M_{+}$ and $M_{-}$. Assume that the following conditions hold
  on the fixed point components $F \subset M^{S^{1}}$:
\begin{align}
\text{If~} F \subset M_{+}, 
  &\quad \Max(\alpha,F) <   w(F); \label{mpluscondition} \\
\text{If~} F \subset M_{-}, 
  &\quad \Min(\alpha,F) >  -w(F). \label{mminuscondition} 
\end{align}
Then 
\[
 Q(M)_{0} = Q(\Mred).
\]
\end{thm}
\begin{proof}
  We will apply Theorem \ref{mainthm} to the manifold with boundary  
  $M_{+}$, whose boundary is $Z$, and apply the ordinary localization
  theorem to all of $M$. Clearly
\[
  Q(\Mred) = \kappa(\alpha|_{M_{+}})
\]
  so Theorem \ref{mainthm} gives
\begin{align}\label{mredcalcwithmainthm}
  Q(\Mred)
    &= \sum_{F \subset M_{+}} (\pi_{F})_{!} \rho 
         \left( \frac{i_{F}^{*}\alpha}{\ee(\nu_{F})} \right) \\
    &= \sum_{F \subset M_{+}} (\pi_{F})_{!} (L_{+,0} - L_{-,0}) 
         \left( \frac{i_{F}^{*}\alpha}{\ee(\nu_{F})} \right) 
\end{align}
while
\begin{align}\label{mzerocalc}
  Q(M)_{0} = ((\pi_{M})_{!} \alpha)_{0}
    &= \left[ \sum_{F \subset M} (\pi_{F})_{!} 
         \left( \frac{i_{F}^{*}\alpha}{\ee(\nu_{F})} \right) \right]_{0} \\
    &= \sum_{F \subset M} (\pi_{F})_{!} (L_{+,0}) 
         \left( \frac{i_{F}^{*}\alpha}{\ee(\nu_{F})} \right) 
\end{align}
so
\begin{equation}\label{compareqrrq}
 Q(\Mred) - Q(M)_{0} 
  = - \sum_{F \subset M_{+}} (\pi_{F})_{!} L_{-,0} 
         \left( \frac{i_{F}^{*}\alpha}{\ee(\nu_{F})} \right) 
    - \sum_{F \subset M_{-}} (\pi_{F})_{!} L_{+,0} 
         \left( \frac{i_{F}^{*}\alpha}{\ee(\nu_{F})} \right). 
\end{equation}
We will analyze the first term; the second is similar.
Recall that by our conventions, the weights of $\nu_{F}$
for $F \subset M_{+}$ are positive. Hence, 
letting $\nu_{F,k}$ be the weight $k$ isotypic
part of $\nu_{F}$, the $\Kso$-Euler class is 
\begin{align*}
 \ee(\nu_{F}) 
  &= \ee(\sum_{k > 0} \nu_{F,k} u^{k}) \\
  &= \prod_{k > 0} \ee(\nu_{F,k} u^{k}) \\
  &= \prod_{k > 0} \sum_{l=0}^{\mathrm{rk}(\nu_{F,k})}
                     (-1)^{l} (\Lambda^{l} \nu_{F,k}) u^{kl}.
\end{align*}
Each factor is a polynomial in $u$ with constant term $1$,
invertible leading coefficient, and degree 
$k \cdot \mathrm{rk}(\nu_{F,k})$. The sum of these degrees is
$w(F)$.

\begin{lemma}\label{invertpolylemma}
Let $R$ be a ring. Let $P \in R[u]$ have 
invertible leading coefficient and constant coefficient and
degree $n$. Then $P$ is invertible in the Laurent ring
$R[[u^{-1}]][u]$, and moreover,
$L_{-}(P^{-1}) \in R[[u^{-1}]][u]$ actually lies in $u^{-n}R[[u^{-1}]]$.
\end{lemma}
\begin{proof}[Proof of the Lemma.]
When $R=\C$, this is clear from examining the poles and order of vanishing  
of $P^{-1}$. The same holds for any ring $R$ by direct calculation:
\begin{align*}
 P &= P_{n}u^{n} + P_{n-1}u^{n-1} + \ldots + P_0 \\  
 P &= P_{n}u^{n} (1 + (P_{n-1}/P_n) u^{-1} + \ldots + (P_0/P_n) u^{-n}) \\  
 P^{-1}    
   &= P_n^{-1} u^{-n} \sum_{l=0}^\infty (-1)^{l}
        ((P_{n-1}/P_n) u^{-1} + \ldots + (P_0/P_n) u^{-n})^l  
        \in u^{-n}R[[u^{-1}]].
\end{align*}
\end{proof}

Hence 
\[
L_{-}(\ee(\nu_{F,k} u^{k}))^{-1} 
   \in u^{k \cdot \mathrm{rk}(\nu_{F,k})}K(F)[[u]]
\]
and 
\[
L_{-}(\ee(\nu_{F}))^{-1} \in u^{w(F)}K(F)[[u]].
\]
Therefore, for $F \subset M_{+}$,
\[
L_{-}\left( \frac{i_{F}^{*}\alpha}{\ee(\nu_{F})} \right)
  \in u^{-w(F)+\Max(\alpha,F)}K(F)[[u^{-1}]].
\]
Given the condition (\ref{mpluscondition}), 
\[
L_{-}\left( \frac{i_{F}^{*}\alpha}{\ee(\nu_{F})} \right)
  \in u^{-1}K(F)[[u^{-1}]]
\]
so 
\[
L_{-,0}\left( \frac{i_{F}^{*}\alpha}{\ee(\nu_{F})} \right) = 0.
\]
Similarly, condition (\ref{mminuscondition}) implies that for
$F \subset M_{-}$,
\[
L_{+,0}\left( \frac{i_{F}^{*}\alpha}{\ee(\nu_{F})} \right) = 0.
\]
Hence 
\[
 Q(\Mred) - Q(M)_{0} = 0.
\]
\end{proof}

This easily implies Theorem \ref{QRthm}.
\begin{proof}[Proof of Theorem \ref{QRthm}.]
$Z = \phi^{-1}(0)$ splits $M$ into $M_{+} = \phi^{-1}((0,\infty))$,
$M_{-} = \phi^{-1}((-\infty,0))$. For any fixed point component
$F$, we have
\[
  i^{*}_{F}\alpha = \alpha_{F} u^{-\phi(F)}
\]
so $\Max(\alpha,F)=\Min(\alpha,F)=-\phi(F)$
which is negative on $M_{+}$, positive on $M_{-}$.
Hence both conditions in Theorem \ref{quantizationtheorem} are 
satisfied, and 
\[
 Q(\Mred) = Q(M)_{0}.
\]
\end{proof}


\end{document}